\documentclass[11pt]{article}
\usepackage[paperwidth=560pt,paperheight=700pt,top=72pt,right=67pt,bottom=72pt,left=72pt]{geometry}
\usepackage{amsmath,amsthm,amsfonts,amssymb,amscd,enumerate,verbatim}
\usepackage[mathscr]{euscript}
\usepackage{array}
\usepackage{hyperref}
\usepackage{fancyhdr}
\usepackage{blindtext}
\usepackage{graphicx,latexsym}
\usepackage{lscape}
\newtheorem{thm}{Theorem}[section]

\newtheorem{defn}[thm]{Definition}

\theoremstyle{definition}

\newtheorem{lemma}{Lemma}[section]

\usepackage{latexsym}
\usepackage{amssymb}
\usepackage{amsmath}
\usepackage{graphicx}
\usepackage{hyperref}
\usepackage{caption}
\usepackage{subcaption}
\usepackage{setspace}

\newcommand{\ind}{1\hspace{-2.3mm}{1}}

\renewcommand{\qed}{\hfill{\ \ \rule{2mm}{2mm}} \vspace{0.2in}}

\begin{document}

\date{}
	\title{Coloured and Dependent Planar Matchings of Random Bipartite Graphs}
		
	\author{\begin{tabular}{ccl}
		\textbf{Ghurumuruhan Ganesan$^{1}$} \\
		\end{tabular}\\
		\begin{tabular}{c}
			$^{1}$IISER, Bhopal\\
			E-mail: gganesan82@gmail.com
	\end{tabular}}
	\maketitle
\pagestyle{fancy}
\lhead{\textit{Planar Matchings}}
\rhead{\thepage}

\vspace{0.2 cm}
\noindent{\bf Abstract:} In this paper, we study two problems related to planar matchings in random bipartite graphs. First, we colour each edge of the complete bipartite graph~\(K_{n,n}\) uniformly randomly from amongst~\(r\) colours and show that if~\(r\) grows linearly with~\(n,\) then the maximum rainbow matching is a non-trivial fraction of~\(r,\) with high probability, i.e. with probability converging to one as~\(n \rightarrow \infty.\) Next we consider planar matchings in a dependent setting where each vertex is forced to choose exactly one neighbour from amongst all possible choices. We obtain estimates for the largest size of a planar matching and also discuss the implication of our results to longest increasing subsequences in enlarged random permutations.\\[0.5em]
{\bf Keywords:} Rainbow Planar Matchings, Dependent Planar Matchings, Random Bipartite Graphs. \\[0.5em]
{\bf{2010 Mathematics Subject Classification:}} Primary:
60J10, 60K35; Secondary: 60C05, 62E10, 90B15, 91D30.






\setcounter{equation}{0}
\renewcommand\theequation{\thesection.\arabic{equation}}
\section{Introduction}
Planar matchings in random graphs have applications in determining longest length of increasing subsequences in permutations and related topics. In~\cite{joh,kiwi,kiwi2} the authors studied  various properties of the largest size of a planar matching in random bipartite graphs with dependencies. In the first part of the paper, we study \emph{rainbow planar matchings} of randomly coloured bipartite graphs. Specifically, we colour each edge of the complete bipartite graph with a random colour and use martingale difference methods and segmentation to obtain deviation bounds on the maximum size of a rainbow planar matching.

In the second part of our paper, we study dependent planar matchings of random bipartite graphs with applications to enlarged random permutations. Random permutations are of great interest from both theoretical and application perspectives. In particular, the longest increasing subsequence~\(M_n\) of a uniformly randomly chosen permutation of~\(\{1,2,\ldots,n\}\) has been well-studied and various properties of~\(M_n\) are known. The paper~\cite{ulam} initiated the study of the longest increasing subsequence of a randomly chosen permutation and~\cite{ham} used subadditive methods to show that the expected value of~\(\frac{M_n}{\sqrt{n}}\) converges to a constant~\(c.\) That~\(c=2\) was independently determined in~\cite{vershik} and~\cite{logan}. Later~\cite{baik} comprehensively determined the asymptotics of~\(M_n\) including a central limit theorem and for a detailed survey, we refer to~\cite{romik}.

In this paper, we consider a variant of uniform permutations which we call as enlarged permutations. Roughly speaking, we increase the ``alphabet size" of the~\(n-\)tuples to~\(k \geq n\) and study the deviation and expectation properties of the longest increasing subsequence~\(T_n = T_n(k).\) The change in the alphabet size affects the underlying distribution due to an increase in size of the sample space and we use a segmentation approach to obtain the desired bounds for~\(T_n.\) Throughout, for completeness, we state and prove our results in the form of planar matchings.

The paper is organized as follows: In Section~\ref{sec_cons_plan}, we state and prove our main result regarding the maximum size of rainbow planar matchings and in Section~\ref{sec_dependent}, we state and prove our main result regarding the maximum size of dependent planar matchings and describe how our result applies to enlarged random permutations.

\setcounter{equation}{0}
\renewcommand\theequation{\thesection.\arabic{equation}}
\section{Rainbow Planar Matchings} \label{sec_cons_plan}
\begin{figure}[tbp]
\centering
\includegraphics[width=3in, trim= 20 450 20 120, clip=true]{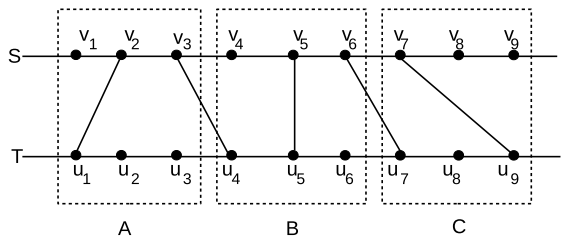}
\caption{Illustration of a planar matching. }
\label{fig_match_ax}
\end{figure}

For~\(i \geq 1,\) let~\(u_i = (i,0)_2\) and~\(v_i = (i,1)_2\) be points in~\(\mathbb{R}^2.\) We use the subscript~\(2\) to differentiate from the two-tuple notations for edges introduced later. For each~\(i \neq j\) join the vertices~\(u_i\) and~\(v_j\) by an edge to obtain an infinite bipartite graph~\(K_{tot}.\)  For integer~\(n \geq 1\) let~\(K_{n,n}\) be the complete bipartite graph containing~\(n\) bottom vertices~\(X = \{u_1,\ldots,u_n\}\) and~\(n\) top vertices~\(Y=\{v_1,\ldots,v_n\}.\) This is illustrated in Figure~\ref{fig_match_ax} where~\(S\) denotes the line~\(x=1\) that contains the vertices of~\(Y\) and~\(T\) denotes the line~\(x=0\) that contains the vertices of~\(X.\) 

An edge~\(e  \in K_{n,n}\) if and only if~\(e\) has one end-vertex~\(u_i \in X\) and the other end-vertex~\(v_j \in Y.\) We denote the edge~\(e\) as~\(e = (u_i,v_j).\) A \emph{matching} of size~\(t\) in~\(K_{n,n}\) is a set of vertex disjoint edges~\({\cal W} = \{e_1,\ldots,e_t\}.\) Suppose~\(e_i\) has end-vertices~\(a_l  = (u_{i_l},1)_2 \in X\) and~\(b_l = (v_{j_l},0)_2 \in Y.\) We say that~\({\cal W}\) is a \emph{planar matching} if~\(i_1<i_2< \ldots <i_t\) and~\(j_1< j_2 < \ldots < j_t.\) In other words, no two edges in~\({\cal W}\) intersect each other. The size of~\({\cal W}\) is defined to the number of edges~\(t.\)

In Figure~\ref{fig_match_ax}, we illustrate the above definition with an example of a planar matching~\({\cal W}_0\) containing six edges~\((u_1,v_2),(u_4,v_3),(u_5,v_5), (u_7,v_6)\) and~\((u_9,v_7).\)

Each edge~\(f\) of~\(K_{n,n}\) is now coloured with a colour~\(X_f\) chosen uniformly randomly from the set~\(\{1,2,\ldots,r\},\) independent of other edges. Let~\(G\) be the resulting randomly coloured graph and let~\({\cal W} = \{h_1,h_2,\ldots,h_t\}\) be a planar matching of~\(G\) containing~\(t\) edges.
\begin{defn}\label{def_2}
We say that~\({\cal W}\) is a \emph{rainbow} planar matching if~\(X_{h_i} \neq X_{h_j}\) for any~\(h_i \neq h_j.\)
\end{defn}
In other words, all the colours in~\({\cal W}\) must be distinct.

Letting~\(R_n\) denote the maximum size of a rainbow planar matching of~\(G,\) we have the following result.
\begin{thm}\label{thm_two} Let~\(r = \alpha n\) for some~\(\alpha > 0.\) There exists~\(0 < \alpha_0 <1\) such that for every~\(\alpha > \alpha_0,\) there are constants~\(0 <\beta_1 \leq \beta_2 < 1\) strictly and~\(\beta_3 > 0\) such that
\begin{equation}\label{dev_est_2}
\mathbb{P}\left(\beta_1 r \leq R_n \leq \beta_2 r\right) \geq 1- e^{-\beta_3 n} \text{ and } \beta_1 r \leq \mathbb{E}R_n \leq \beta_2 r.
\end{equation}
Moreover,~\(var(R_n) \leq 2 \mathbb{E}R_n.\)
\end{thm}
In words, the above result says that the largest size of the rainbow planar matching is a non-trivial fraction of the total number of colours~\(r,\) if the parameter~\(\alpha\) is close to~\(1.\)  On the other hand, the proof of our result can be also adapted to show that if~\(r = o(n),\) then~\(\frac{R_n}{r} \longrightarrow 1\) in probability as~\(n \rightarrow \infty.\) This hints at the  possibility of a transition point for the size of the rainbow planar matching, with respect to the parameter~\(\alpha.\)






\emph{Proof Theorem~\ref{thm_two}}:  The proof of the variance bound uses a pivotal edge argument and is exactly as in the strong matching case (see proof of Theorem~\(2.2,\)~\cite{gan}). Next we show that the expectation bounds follow from the deviation bounds. Indeed from the deviation lower bound in~(\ref{dev_est_2}), we get~\(\mathbb{E}R_n \geq \beta_1 r (1-e^{-\beta_3 n})\) and using the fact that~\(R_n \leq r\) and the upper bound in~(\ref{dev_est_2}), we get that~\(\mathbb{E}R_n \leq \beta_2 r + r e^{-\beta_3 n}.\) This obtains the expectation bounds in~(\ref{dev_est_2}).

We now begin with the proof of the upper deviation bound in~(\ref{dev_est_2}). Let~\({\cal S} = \{(u_i,v_i)\}_{1 \leq i \leq t}\) be any planar matching containing~\(t\) edges. The probability that~\({\cal S}\) is a rainbow planar matching is~\(\prod_{i=1}^{t-1}\left(\frac{r-i}{r}\right) \) which is bounded above by
\begin{eqnarray}
\exp\left(-\frac{1}{r}\sum_{i=1}^{t-1}i \right) &=& \exp\left(-\frac{t(t-1)}{2r}\right) \nonumber\\
&=& \exp\left(\frac{t}{2r}\right)\exp\left(-\frac{t^2}{2r}\right) \nonumber\\
&\leq& \sqrt{e}\exp\left(-\frac{t^2}{2r}\right), \nonumber
\end{eqnarray}
since~\(t \leq r.\)  The number of choices for~\(\{u_i\}_{1 \leq i \leq t}\) is~\({n \choose t}\) and an analogous bound holds for~\(\{v_i\}_{1 \leq i \leq t}.\) Therefore if~\(F_n = F_n(\epsilon), \epsilon < \frac{1}{2}\) denotes the event that every rainbow planar matching has size at most~\((1-\epsilon)r,\) then by the union bound we have that
\begin{equation}\label{fn_c}
\mathbb{P}(F_n^c) \leq \sqrt{e}\sum_{t \geq (1-\epsilon)r+1} {n \choose t}^2 \cdot \exp\left(-\frac{t^2}{2r}\right) \leq \sqrt{e}\sum_{t \geq (1-\epsilon)r+1} {n \choose t}^2 e^{-(1-\epsilon)^2r/2}.
\end{equation}

Assuming~\(\alpha > \frac{1}{2},\) we choose~\(\epsilon > 0\) small enough to apply the monotonicity of the Binomial coefficient and the Stirling approximation successively to get that~\({n \choose t} \leq {n \choose (1-\epsilon)r} \leq ne^{nH((1-\epsilon)\alpha)}\)
where~\(H(x) = -x \log{x} - (1-x) \log(1-x)\) is the binary entropy function and the logarithms are natural throughout. From~(\ref{fn_c}) we therefore get that
\begin{equation}\label{fnc_ax}
\mathbb{P}(F_n^c) \leq  \sqrt{e}\sum_{t \geq (1-\epsilon)r+1} n^2e^{2nH((1-\epsilon)\alpha)}e^{-(1-\epsilon)^2\alpha n/2}.
\end{equation}
We know that~\(H(x) \longrightarrow 0\) as~\(x \rightarrow 0.\) Therefore if~\(\alpha_0\) is the solution to the equation~\(2H(x) = \frac{x}{2},\) then for every~\(\alpha > \alpha_0,\) we can choose~\(\epsilon > 0\) small in~(\ref{fnc_ax}) to get that~\(\mathbb{P}(F_n^c) \leq e^{-\delta n}\)
for some constant~\(\delta  >0.\) Fixing such an~\(\epsilon,\) gives the upper deviation bound in~(\ref{dev_est_2}).

For the lower deviation bound in~(\ref{dev_est_2}), we argue as follows: For~\(0 < \epsilon < \frac{1}{2},\) let~\(E_n\) be the event that the maximum rainbow planar matching contains~\(t \leq \epsilon r\) edges. We bound~\(\mathbb{P}(E_n)\) by estimating the size of the colour set of the~\(n\) vertical edges~\(f_i = (u_i,v_i), 1 \leq i \leq n.\)  The probability that the colour~\(X_{f_i} \) of the edge~\(f_i\) belongs to the set~\(\{1,2,\ldots,t\}\) equals~\(\frac{t}{r} \leq \epsilon.\) Therefore the probability that each~\(f_i,1 \leq i \leq n\) is coloured with a colour from~\(\{1,2,\ldots,t\}\) is at most~\(\epsilon^{n}.\) Since there are~\({r \choose t} \leq {r \choose \epsilon r} \leq 2^{r} =2^{\alpha n}\) ways of choosing~\(t\) colours from the set of all available colours, we get by the union bound that
\[\mathbb{P}(E_n) \leq \sum_{1 \leq t \leq \epsilon r} 2^{\alpha n}\epsilon^{n} \leq n (2^{\alpha} \epsilon)^{n} \leq e^{-\gamma n}\]
for some constant~\(\gamma > 0\) provided~\(\epsilon > 0\) is small enough. Fixing such an~\(\epsilon,\) we get the lower deviation bound in~(\ref{dev_est_2}).~\(\qed\)

\setcounter{equation}{0}
\renewcommand\theequation{\thesection.\arabic{equation}}
\section{Dependent Planar Matchings}\label{sec_dependent}
Let~\(K_{tot}\) be the infinite bipartite graph described in Section~\ref{sec_cons_plan} and for integers~\(k \geq n \geq 1,\) let~\(K_{k,n} \subseteq K_{tot}\) be the complete bipartite graph  containing the~\(k\) top vertices~\(v_1,\ldots,v_k\) and the~\(n\) bottom vertices~\(u_1,\ldots,u_n.\) Let~\({\cal G}_{k,n}\) be the set of all~\(1-\)regular subgraphs of~\(K_{k,n};\) i.e., the set of all subgraphs of~\(K_{k,n}\) satisfying the property that each bottom vertex~\(u_i, 1 \leq  i \leq n\) is adjacent to exactly one top vertex~\(v_j, 1 \leq j \leq k\) and each top vertex is adjacent to at most one bottom vertex. Let~\(G\) be uniformly randomly chosen from~\({\cal G}_{k,n}\) and let~\(T_n = T_n(k)\) be the largest size of a planar matching of~\(G.\)

For~\(k = n,\) the quantity~\(T_n\) could be interpreted as the length of the longest increasing subsequence in a uniformly randomly chosen permutation of~\(n.\) In this case, it is well-known that~\(T_n\) is of the order of~\(\sqrt{n}\) with high probability~\cite{romik}, i.e. with probability converging to one as~\(n \rightarrow \infty.\) We could therefore think of the above dependent planar matching setup as an ``enlarged permutation" where we choose~\(n-\)tuples with distinct entries from the set~\(\{1,2,\ldots,k\}\) with~\(k \geq n.\)

For general~\(k \geq n,\) we have the following result regarding the mean and variance of~\(T_n.\) As before, constants do not depend on~\(n.\)
\begin{thm}\label{thm_dep} We have:\\
\((a)\) For every~\(\epsilon  >0,\) there are positive constants~\(C_i, 0 \leq i \leq 3\) and a constant~\(0 < D_0 < 1\) such that for all~\(n \geq C_0,\) \begin{equation}\label{mean_bounds}
\mu_{low} := (1-e^{-1}) \sqrt{n} - \frac{17}{\sqrt{n}} \leq \mathbb{E}T_n \leq (e+\epsilon)\sqrt{n} + 1 =: \mu_{up},
\end{equation}
\begin{equation}\label{dev_bounds}
\mathbb{P}\left(T_n \leq \mu_{up}\right) \geq 1-e^{-C_1 \sqrt{n}} \text{ and } \mathbb{P}\left(T_n \geq (1-\epsilon) \mu_{low}\right) \geq 1-z_n,
\end{equation}
where~\(z_n:=\min\left(D_0,\frac{C_2}{\sqrt{n}}+\frac{C_2n}{k}\right) .\)\\
\((b)\) For every~\(2< b = b(n) \leq \frac{n}{(32)^2\log{n}}\) we have that
\begin{equation}\label{low_bounds_gen}
\mathbb{P}\left(T_n \geq \sqrt{\frac{n}{b\log{n}}}\right) \geq 1-\exp\left(-\left(\frac{b}{2}-1\right)\log{n}\right).
\end{equation}
\end{thm}
From~(\ref{dev_bounds}) we see that if~\(k\) grows faster than~\(n,\) then the maximum planar matching size~\(T_n\) is of the order of~\(\sqrt{n}\) with high probability, i.e. with probability converging to one as~\(n \rightarrow \infty.\)






The proof of the upper bounds in Theorem~\ref{thm_dep} is analogous to the case of random permutations (see Lemma~\(1.5\) in Romik (2014)). For completeness, we give small proofs in the next section.  In our proof of the lower deviation bounds in Theorem~\ref{thm_dep}, we use segmentation to split the given graph~\(K_{k,n}\) into smaller subgraphs and use planar matchings within the smaller subgraphs to obtain an estimate on the overall quantity~\(T_n.\) Therefore, in the next subsection, we collect preliminary results regarding the minimum size of a segmented planar matching of~\(K_{k,n}\) and independent planar matchings, respectively. In the final subsection, we prove Theorem~\ref{thm_dep}. 

\subsection*{Segmented Planar Matchings}
For integer~\(1 \leq t \leq n\) and for~\(1 \leq i \leq I := \frac{n}{t}\) let~\({\cal T}_i\) be the set of all integers~\(l\) satisfying~\((i-1)t+1 \leq l \leq i t.\) For simplicity we have assumed that~\(\frac{n}{t}\) is an integer; if not then we let~\(I\) be the largest integer such that~\(I \cdot t \leq n\) and allow the last ``segment"~\({\cal T}_I\) to be the set of all integers~\(j\) satisfying\\\((I-1) t+1 \leq i \leq n.\) In this case,~\({\cal T}_I\) contains between~\(\frac{n}{t}\) and~\(\frac{2n}{t}\) integers. Similarly for an integer~\(1 \leq s \leq k\) and for~\(1 \leq j \leq J := \frac{k}{s}\) let~\({\cal S}_i\) be the set of all integers~\(l\) satisfying~\((j-1)s+1 \leq l \leq j s.\)  As before, we assume that~\(\frac{k}{s}\)  is an integer and we always choose~\(s \leq \frac{tk}{n}\) so that the number of ``bottom" segments~\(I \leq J,\) the number of ``top" segments.

Suppose~\(s = \frac{tk}{n}\) so that~\(I = J\) and let~\({\cal W}\) be a planar matching of~\(G.\) We say that~\({\cal W}\) is a~\(t-\)segmented planar matching if for each edge~\(e \in {\cal W},\) there exists an integer~\(1 \leq i \leq I\) such that one end-vertex of~\(e\) belongs to~\({\cal T}_i\) and the other end-vertex of~\(e\) belongs to~\({\cal S}_i.\) For~\(1 \leq i \leq \frac{n}{t},\) let~\(A_i\) be the event that there exists an integer~\(z \in {\cal T}_i\) such that~\(\pi(z) \in {\cal S}_i.\) Defining the \emph{minimum size} of a~\(t-\)segmented planar matching of the random graph~\(G\) to be
\begin{equation}\label{x_def}
X_t  = X_t(s) := \sum_{i=1}^{n/t} \ind(A_i),
\end{equation}
we have the following estimates for the mean and variance of~\(X_t.\)
\begin{lemma}\label{thm_seg} Suppose~\(s = \frac{kt}{n}\) so that~\(I = J.\) Letting~\(\mu_t := \frac{n}{t} \left(1-\exp\left(-\frac{t^2}{n}\right) \right),\) we have that
\begin{equation}\label{mean_xt}
\mu_t - \frac{17}{t}e^{-\frac{t^2}{n}} \leq \mathbb{E}X_t \leq \mu_t + \frac{32n}{k^2}e^{-\frac{t^2}{n}}.
\end{equation}
If~\(\theta_1 \sqrt{n} \leq t \leq \theta_2 \sqrt{n}\) for some constants~\(\theta_1,\theta_2 > 0,\) then there exists a positive constant~\(D = D(\theta_1,\theta_2)\) such that
\begin{equation}\label{var_xt}
var(X_t) \leq D \left(\sqrt{n} + \frac{n^2}{k}\right)
\end{equation}
\end{lemma}
In the proof of Theorem~\ref{thm_dep} in the next Section, we use the fact that~\(X_t\) is a lower bound on~\(T_n\) to obtain deviation and expectation estimates for~\(T_n.\)


\emph{Proof of~(\ref{mean_xt}) in Lemma~\ref{thm_seg}}: We use Lemma~\ref{pa_lem} in Appendix with~\(s= \frac{kt}{n}\) to get from the upper bound in~(\ref{pa_one}) that
\[\mathbb{P}(A_1^c) \leq \exp\left(-\frac{t^2}{n} + \frac{4(s+t)^2t}{k^2}\right).\]
The term
\begin{equation}\label{gen_st}
\frac{(s+t)^2t}{k^2} = \frac{(k+n)^2t}{k^2n^2} \leq \frac{4t}{n^2} \leq \frac{4}{n}
\end{equation}
where the first inequality in~(\ref{gen_st}) is true since~\(k \geq n\) and the second inequality in~(\ref{gen_st}) follows from the fact that~\(t \leq n.\) Thus
\begin{equation}\label{pa_one_estax}
\mathbb{P}(A_1^c) \leq \exp\left(-\frac{t^2}{n} + \frac{16}{n}\right) \leq \exp\left(-\frac{t^2}{n}\right) \left(1+\frac{17}{n}\right)
\end{equation}
for all~\(n\) large and so
\[\mathbb{E}X_t = \frac{n}{t} \left(1-\mathbb{P}(A_1^c)\right)\geq \frac{n}{t}\left(1-\exp\left(-\frac{t^2}{n}\right)\right) - \frac{17}{t}e^{-\frac{t^2}{n}}.\] This obtains the lower bound in~(\ref{mean_xt}).

Similarly using the lower bound in~(\ref{pa_one}) we get that
\[\mathbb{P}(A_1^c) \geq \exp\left(-\frac{t^2}{n} - \frac{8(s+t)^2t}{k^2}\right),\]
where~\(\frac{(s+t)^2t}{k^2} = \frac{(k+n)^2t}{k^2n^2} \geq \frac{4t}{k^2}\) since~\( k \geq n.\) Thus
\[\mathbb{P}(A_1^c) \geq \exp\left(-\frac{t^2}{n}\right) \exp\left( -\frac{32t}{k^2}\right) \geq \exp\left(-\frac{t^2}{n}\right) \left( 1- \frac{32t}{k^2}\right)\] and consequently,
\[\mathbb{E}X_t = \frac{n}{t} \left(1-\mathbb{P}(A_1^c)\right) \leq \frac{n}{t} \left(1-\exp\left(-\frac{t^2}{n}\right) \right) + \frac{32 n}{k^2} e^{-\frac{t^2}{n}}.\]
This obtains the upper bound in~(\ref{mean_xt}).~\(\qed\)


\emph{Proof of~(\ref{var_xt}) in Lemma~\ref{thm_seg}}: To obtain the deviation estimate for~\(X_t,\) we let~\(Y = \frac{n}{t}-X_t\) and first estimate the variance of~\(Y.\) Setting
\[cov(\ind(A_i^c),\ind(A_j^c))  := \mathbb{P}(A_i^c \cap A_j^c) - \mathbb{P}(A_i^c) \mathbb{P}(A_j^c)  \]
to be the covariance between the indicator functions of the events~\(A_i\) and~\(A_j,\) we have that
\begin{eqnarray}
var(Y) &=& \sum_{i=1}^{n/t} \mathbb{P}(A^c_i) - \mathbb{P}^2(A^c_i) + 2\sum_{i <j} cov(\ind(A_i^c),\ind(A_j^c)) \nonumber\\
&\leq& \frac{n}{t}\left(\mathbb{P}(A^c_1) - \mathbb{P}^2(A^c_1)\right) + \frac{2n^2}{t^2} cov(\ind(A_1^c),\ind(A_2^c))\nonumber\\
&\leq& \frac{n}{t}\mathbb{P}(A^c_1)  + \frac{2n^2}{t^2} cov(\ind(A_1^c),\ind(A_2^c)) \nonumber\\
&\leq& \frac{2n}{t}e^{-\frac{t^2}{n}}+ \frac{2n^2}{t^2} cov(\ind(A_1^c),\ind(A_2^c)) \label{var_x_est}
\end{eqnarray}
for all~\(n\) large, using~(\ref{pa_one_estax}).

For the second term in~(\ref{var_x_est}), we use the correlation estimate~(\ref{pa_two}) in Lemma~\ref{pa_lem} to get that
\[cov(\ind(A_1^c),\ind(A_2^c)) \leq \mathbb{P}(A_1^c) \mathbb{P}(A_2^c)\left(e^{\frac{5t^2}{k}}-1\right).\]
Since~\(t \leq \theta_2 \sqrt{n},\) we have that~\(\frac{5t^2}{k} \leq \frac{5\theta_2^2n}{k} \leq 5\theta_2^2.\)
For all~\(x \leq 5\theta_2^2,\) we have that~\(e^{x} \leq 1+D_1x\) for some large constant~\(D_1>0\) and so again using~\(\theta_1 \sqrt{n} \leq t \leq \theta_2 \sqrt{n},\) we get that
\begin{equation}\label{cov_est}
cov(\ind(A_1^c),\ind(A_2^c))  \leq 5D_1 \frac{t^2}{k} \mathbb{P}(A_1^c) \mathbb{P}(A_2^c)  \leq \frac{D_2 n}{k} e^{-\frac{2t^2}{n}} \leq \frac{D_3 n}{k}
\end{equation}
for some constants~\(D_2,D_3 > 0.\) Plugging~(\ref{cov_est}) into~(\ref{var_x_est}) and using~\(\theta_1 \sqrt{n} \leq t \leq \theta_2 \sqrt{n},\) we get that
\begin{equation}\nonumber
var(Y) \leq \frac{2n}{t}e^{-\frac{t^2}{n}} + \frac{n^2}{t^2} \cdot\frac{D_3n}{k} \leq D_4 \left(\sqrt{n} + \frac{n^2}{k}\right)
\end{equation}
for some constant~\(D_4 > 0.\) Since~\(var(X_t) = var(Y),\) this completes the proof of~(\ref{var_xt}).~\(\qed\)

\subsection*{Proof of Theorem~\ref{thm_dep}}
\emph{Proof of the upper bound in~(\ref{mean_bounds})}: For a given integer~\(l \geq 1\) we first estimate the probability of the event~\(E_l\) that the random graph~\(G\) contains a planar matching of size~\(l.\) Indeed, let~\(\pi(i)\) be the~\(x-\)coordinate of the top vertex adjacent to the bottom vertex~\(u_i\) in the graph~\(G\) and suppose that~\(\pi({i_1}) < \ldots < \pi(i_l)\) for some deterministic integers~\(1 \leq i_1 < i_2 < \ldots <i_l \leq n.\) The number of choices for~\(\{i_1,\ldots,i_l\}\) is~\({ n \choose l}\) and the number of choices for~\(\{\pi({i_1}),\ldots,\pi({i_l})\}\) is~\({k \choose l}.\) From~(\ref{pi_basic}) and the union bound, we then get that
\[\mathbb{P}(E_l) \leq {n \choose l} \cdot {k \choose l} \cdot \frac{1}{k(k-1)\ldots (k-l+1)} = {n \choose l} \cdot \frac{1}{l!}.\]
Using~\({n \choose l} \leq \left(\frac{ne}{l}\right)^{l}\) and~\(l! \geq l^{l}e^{-l},\) we get that~\(\mathbb{P}(E_l) \leq \left(\frac{ne^2}{l^2}\right)^{l} \leq \left(\frac{e}{a}\right)^{2l}\) for~\(l \geq a \sqrt{n}.\)

If~\(a > e\) strictly then~\(\left(\frac{e}{a}\right)^{2l} =e^{-2\delta l}\) for some constant~\(\delta > 0\) and so\\\( \mathbb{P}(E_l) \leq e^{-2\delta l}.\) If~\(E_{up}\) denotes the event that~\(G\) contains a planar matching of length at least~\(a \sqrt{n},\) then from the union bound, we get that
\begin{equation}\label{e_up_est}
\mathbb{P}(E_{up}) \leq n \cdot e^{-2a \delta \sqrt{n}} \leq e^{-a \delta \sqrt{n}}
\end{equation}
for all~\(n\) large. To upper bound the expectation, we use the fact that if~\(E_{up}\) occurs, then~\(T_n \leq a\sqrt{n}\) and if the complement event~\(E_{up}^c\) occurs, then~\(T_n \leq n.\) Therefore from~(\ref{e_up_est}), we get that~\(\mathbb{E}T_n \leq a\sqrt{n} + n \cdot e^{-a\delta \sqrt{n}} \leq a\sqrt{n}+1\) for all~\(n\) large. This obtains the upper bound in~(\ref{mean_bounds}). ~\(\qed\)

\emph{Proof of the lower bound in~(\ref{mean_bounds})}: By definition the term~\(X_t\) defined in~(\ref{x_def}) is a lower bound on the maximum size~\(T_n\) of a planar matching of~\(G.\) To see this is true suppose~\(X = w\) and suppose~\(A_{j_1},\ldots,A_{j_w}\) occur with~\(j_1 < \ldots < j_w.\) For each~\(1 \leq u \leq w,\) there exists an index~\(p_u \in {\cal T}_{j_u}\) and an integer~\(q_u \in {\cal S}_{j_u}\) such that~\(\pi(p_u) = q_u.\) By definition~\(\pi(p_1) < \pi(p_2) < \ldots < \pi(p_w)\) and so~\(G\) contains a planar matching of size at least~\(w.\)

We therefore use Lemma~\ref{thm_seg} with~\(t = \sqrt{n}\) to get that~\(\mu_t = \sqrt{n}(1-e^{-1})\) and so from the lower bound in~(\ref{mean_xt}), we get that~\(\mathbb{E}X_t \geq (1-e^{-1})\sqrt{n} -\frac{17}{\sqrt{n}}.\) This completes the proof of the lower bound in~(\ref{mean_bounds}).~\(\qed\)



\emph{Proof of~(\ref{dev_bounds})}: The upper bound follows from~(\ref{e_up_est}). To prove the lower bound, we use Lemma~\ref{thm_seg} with~\(t = \sqrt{n}.\) First from the upper bound in~(\ref{mean_xt}) and the fact that~\(k \geq n,\) we get that~\[\mathbb{E}X_t \leq (1-e^{-1})\sqrt{n} + \frac{32n}{k^2} \leq (1-e^{-1})\sqrt{n}+ \frac{32}{n}\leq \sqrt{n}.\]
From the variance estimate~(\ref{var_xt}), we already know that
\begin{equation}\label{var_est_ax_dep}
var(X_t) \leq D \left(\sqrt{n} + \frac{n^2}{k}\right).
\end{equation}

We now use the Paley-Zygmund inequality to show that~\(X_t \geq (1-\epsilon)\mu_{low}\) with positive probability and then use the fact that~\(X_t\) is a lower bound for~\(T_n\) to get that~\(T_n \geq (1-\epsilon) \mu_{low}\) with positive probability.  From Paley-Zygmund inequality we have for~\(\epsilon > 0\) that
\[\mathbb{P}(X_t \geq (1-\epsilon) \mathbb{E}X_t) \geq \frac{\epsilon^2 (\mathbb{E}X_t)^2}{\mathbb{E}X^2_t}\] and so using~\(\mathbb{E}X_t \geq \mu_{low} \geq \frac{\sqrt{n}}{2},\) we have~\(\mathbb{P}(X_t \geq (1-\epsilon) \mu_{low}) \geq \frac{\epsilon^2n}{4\mathbb{E}X^2_t}.\) From~(\ref{var_est_ax_dep}) and the fact that~\(k \geq n,\) we then obtain~\(\mathbb{E}X^2_t \leq D_1 n\) for some constant~\(D_1 > 0.\)  Consequently~\(X_t \geq (1-\epsilon)\mu_{low}\) with probability at least~\(\frac{\epsilon^2}{4D_1} =: D_0.\)

Next, using the Chebychev's inequality, we have for~\(\epsilon > 0\) that
\[\mathbb{P}\left(|X_t-\mathbb{E}X_t| \leq \epsilon \mathbb{E}X_t\right) \geq 1-\frac{var(X_t)}{\epsilon^2 (\mathbb{E}X_t)^2}.\]
Since~\(\mathbb{E}X_t \geq \mu_{low}\) we get from~(\ref{var_est_ax_dep}) that~\(\mathbb{P}\left(X_t \geq (1-\epsilon)\mu_{low}\right) \geq 1-\frac{D_2}{\sqrt{n}} - \frac{D_3n}{k}\) for some positive constants~\(D_i = D_i(\epsilon), i=2,3.\) This obtains the second lower bound in~(\ref{dev_bounds}) and therefore completes the proof of the deviation lower bound in~(\ref{dev_bounds}).~\(\qed\)

\emph{Proof of~(\ref{low_bounds_gen})}: We set~\(t := \sqrt{bn \log{n}}\) and~\(s := \frac{kt}{n}\) so that the number of segments~\(I=J\) and use Lemma~\ref{pa_lem} to estimate the probability of the event~\(\mathbb{P}(A_i) = \mathbb{P}(A_1).\) By the definition of~\(t\) and the fact that~\(k \geq n,\) we respectively get
\begin{equation}\label{st_eval}
\frac{st}{k} = \frac{t^2}{n} = b\log{n} \text{ and } \frac{s+t}{k} = \frac{t}{k} + \frac{t}{n} \leq \frac{2t}{n} = o(1)
\end{equation}
and so the conditions in Lemma~\ref{pa_lem} are satisfied. Moreover, using~\(k \geq n\) again, we also get that
\begin{equation}\label{rem_eval}
\frac{(s+t)^2 t}{k^2} = \frac{t^3(k+n)^2}{k^2n^2} \leq \frac{4t^3}{n^2} = \frac{4(b\log{n})^{3/2}}{\sqrt{n}}
\end{equation}

Plugging~(\ref{st_eval}) and~(\ref{rem_eval}) into the upper bound of~(\ref{pa_one}), we get that
\begin{eqnarray}\label{eva}
\mathbb{P}(A_1^c) &\leq& e^{-b\log{n}} \cdot \exp\left(\frac{16(b\log{n})^{3/2}}{\sqrt{n}}\right) \nonumber\\
&=& \exp\left(-b\log{n}\left(1-16\sqrt{\frac{b\log{n}}{n}}\right)\right) \nonumber\\
&\leq& \exp\left(-\frac{b\log{n}}{2}\right) \nonumber
\end{eqnarray}
for all~\(n\) large since~\( b \leq \frac{n}{(32)^2\log{n}}\) (see statement of Theorem~\ref{thm_dep}) and so by the union bound we get that~\(\bigcap_{i=1}^{n/t}A_i\) occurs with probability at least~\[1-\frac{n}{t} \cdot \exp\left(-\frac{b}{2} \log{n}\right) \geq 1-n \cdot \exp\left(-\frac{b}{2} \log{n}\right).\] This in turn implies that~\(T_n \geq X_t = \frac{n}{t} = \sqrt{\frac{n}{b\log{n}}}\) with probability least\\\(1-n \cdot \exp\left(-\frac{b}{2} \log{n}\right),\) completing the proof of~(\ref{low_bounds_gen}).~\(\qed\)




\setcounter{equation}{0}
\renewcommand\theequation{A.\arabic{equation}}
\section{Conclusion} \label{conc}
In this paper, we have studied coloured and dependent planar matchings in random bipartite graphs. In the colouring part, we have shown that the largest rainbow matching is a non-trivial fraction of the total number of colours with high probability. This hints at a possible transition point with respect to the colouring parameter. In the dependent setting, we have obtained estimates on the largest size of a planar matching and also explained our results in terms of longest increasing subsequences in random enlarged permutations.

\setcounter{equation}{0}
\renewcommand\theequation{A.\arabic{equation}}
\section*{Appendix}
Throughout we use the following standard deviation estimate.
\begin{lemma} \label{lem_conc}
Let~\(\{X_j\}_{1 \leq j \leq r}\) be independent Bernoulli random variables with~\[\mathbb{P}(X_j = 1) = 1-\mathbb{P}(X_j = 0) > 0.\] If~\(T_r := \sum_{j=1}^{r} X_j,\theta_r := \mathbb{E}T_r\) and~\(0 < \gamma \leq \frac{1}{2},\) then
\begin{equation}\label{conc_est_f}
\mathbb{P}\left(\left|T_r - \theta_r\right| \geq \theta_r \gamma \right) \leq 2\exp\left(-\frac{\gamma^2}{4}\theta_r\right)
\end{equation}
for all \(r \geq 1.\)
\end{lemma}
For a proof of~(\ref{conc_est_f}), we refer to Corollary A.1.14, pp. 312 of Alon and Spencer (2008).


To prove Lemma~\ref{thm_seg}, we use the following preliminary result that obtains estimates on the probability of the events~\(A_i^c\) and~\(A_i^c \cap A_j^c.\) 
\begin{lemma}\label{pa_lem} Let~\(1 \leq t \leq n\) and~\(1 \leq s \leq k\) be any two integers satisfying\\\(\frac{s+t}{k} \leq \frac{1}{8}.\) We have that
\begin{equation}\label{pa_one}
\exp\left(-\frac{st}{k}-\frac{8(s+t)^2t}{k^2}\right) \leq \mathbb{P}(A_1^c)  \leq \exp\left(-\frac{st}{k}+\frac{4(s+t)^2t}{k^2}\right)
\end{equation}
and
\begin{equation}\label{pa_two}
\mathbb{P}(A_1^c \cap A_2^c) \leq \mathbb{P}(A_1^c) \mathbb{P}(A_2^c) \exp\left(\frac{5t^2}{k}\right).
\end{equation}
\end{lemma}

\emph{Proof of~(\ref{pa_one}) in Lemma~\ref{pa_lem}}: For any integer~\(1 \leq l \leq n\) and~\(l\) distinct integers~\(1 \leq a_1,\ldots,a_l \leq n,\) we have that
\begin{equation}\label{pi_basic}
\mathbb{P}(u_1 = a_1,\ldots,u_l=a_l) = \frac{1}{k(k-1)\ldots(k-l+1)}.
\end{equation}
Let~\(\pi(i)\) be the~\(x-\)coordinate of the top vertex adjacent to~\(u_i\) in the graph~\(G.\) If the event~\(A_1^c\) occurs, then the number of choices for~\(\pi(1)\) is~\(k-s.\) Similarly, given~\(\pi(1),\) the number of choices for~\(\pi(2)\) is~\(k-s-1\) and so on. Thus the number of choices for the~\(t-\)tuple~\((\pi(1),\ldots,\pi(t))\) is~\[(k-s)\cdot (k-s-1)\cdots(k-s-t+1)\] and so we get from~(\ref{pi_basic}) that
\begin{eqnarray}
\mathbb{P}(A_1^c) &=& \frac{(k-s) \cdot (k-s+1) \ldots (k-s-t+1)}{k(k-1)\cdots(k-t+1)} \nonumber\\
&=& \frac{\left(1-\frac{s}{k}\right) \cdots \left(1-\frac{s+t-1}{k}\right)}{\left(1-\frac{1}{k}\right) \cdots \left(1-\frac{t-1}{k}\right)}. \label{eq_one}
\end{eqnarray}



Using~\(\frac{s+t}{k} \leq \frac{1}{4}\) we get for all~\(x \leq \frac{2(s+t)}{k}\) that
\begin{eqnarray}
\frac{e^{-x}}{1-x} &=& 1+\frac{1}{1-x} \left(\frac{x^2}{2!}-\frac{x^3}{3!} + \ldots \right) \nonumber\\
&\leq& 1 + \frac{2x^2}{1-x} \nonumber\\
&\leq& 1+ 4x^2.\label{karla}
\end{eqnarray}
Thus
\begin{equation}\label{ex_est}
e^{-x} \cdot \left(1+\frac{4(s+t)^2}{k^2}\right)^{-1} \leq 1-x \leq e^{-x}.
\end{equation}

Plugging the upper bound of~(\ref{ex_est}) into the numerator of~(\ref{eq_one}) and the lower bound of~(\ref{ex_est}) into the denominator, we get that~\(\mathbb{P}(A_1^c) \leq I_1 \cdot I_2\) where
\begin{equation}\label{eq_two}
I_1 := \frac{\exp\left(-\frac{1}{k}\sum_{l=s}^{s+t-1}k\right)}{\exp\left(-\frac{1}{k}\sum_{l=1}^{t-1}k\right)}  = \exp\left(-\frac{st}{k}\right)
\end{equation}
and
\[I_2 := \left(1+\frac{4(s+t)^2}{k^2}\right)^{t} \leq \exp\left(\frac{4(s+t)^2t}{k^2}\right).\] This obtains the upper bound in~(\ref{pa_one}).

Similarly, substituting the lower bound of~(\ref{ex_est}) into the numerator of~(\ref{eq_one}) and the upper bound of~(\ref{ex_est}) into the denominator, we get that
\begin{equation}\label{temp_karla}
\mathbb{P}(A_1^c) \geq e^{-\frac{st}{k}} \cdot \left(1+\frac{4(s+t)^2}{k^2}\right)^{-t} \geq e^{-\frac{st}{k}}\cdot \left(1-\frac{4(s+t)^2}{k^2}\right)^{t}.
\end{equation}
Again using~\(\frac{s+t}{k} \leq \frac{1}{2}\) and the fact that~\(1-x\geq e^{-2x}\) for all~\(x< \frac{1}{2},\) we get the lower bound in~(\ref{pa_one}) from the final expression in~(\ref{temp_karla}).~\(\qed\)

\emph{Proof of~(\ref{pa_two}) in Lemma~\ref{pa_lem}}:  We use a split set argument and as before,~\(\pi(i), 1 \leq i \leq n\) is the~\(x-\)coordinate of the top vertex adjacent to~\(u_i,\) in the graph~\(K_{k,n}.\) For a deterministic set~\({\cal Q}_1 \subseteq \{1,2,\ldots,t\},\) let~\(F_1({\cal Q}_1)\) be the event that~\(\pi(j) \in \{s+1,\ldots,2s\}\) for every integer~\(j \in {\cal Q}_1\) and~\(\pi(l) \notin \{1,2,\ldots,2s\}\) for every~\(l \in \{1,2,\ldots,t\} \setminus {\cal Q}_1.\) Similarly for~\({\cal Q}_2 \subseteq \{t+1,\ldots,2t\},\) let~\(F_2({\cal Q}_2)\) be the event that~\(\pi(j) \in \{1,\ldots,s\}\) for every integer~\(j \in {\cal Q}_2\) and~\(\pi(j) \notin \{1,2,\ldots,2s\}\) for every~\(j \in \{t+1,\ldots,2t\} \setminus {\cal Q}_2.\) By definition we then have that
\begin{equation}\label{ea_split}
A_1^c \cap A_2^c = \bigcup_{{\cal Q}_1,{\cal Q}_2} F_1({\cal Q}_1) \cap F_2({\cal Q}_2),
\end{equation}
where the union is over all sets~\({\cal Q}_1 \subseteq \{1,2,\ldots,t\}\) and~\({\cal Q}_2 \subseteq \{t+1,\ldots,2t\}.\) Thus
\begin{equation}\label{ea_split_ax}
\mathbb{P}(A_1^c \cap A_2^c) = \sum_{j_1,j_2}\;\;\;\sum_{\#{\cal Q}_1 = j_1} \sum_{\#{\cal Q}_2 = j_2} \mathbb{P}\left(F_1({\cal Q}_1) \cap F_2({\cal Q}_2)\right).
\end{equation}

Letting~\((a)_b :=
\left\{
\begin{array}{cc}
a(a-1)\cdots(a-b+1) & \text{ if }b  \geq 1\\
1 & \text{ if } b = 0,
\end{array}
\right.\)
we have that~\[\mathbb{P}(F_1({\cal Q}_1)) = \frac{(s)_{j_1}(k-2s)_{t-j_1}}{(k)_t}.\] A similar expression holds for~\(F_2({\cal Q}_2)\) and the intersection~\(F_1({\cal Q}_1) \cap F_2({\cal Q}_2).\) Letting~\(j_{low} := \min(j_1,j_2)\) and~\(j_{up} := \max(j_1,j_2)\) and combining the above expressions, we get that
\begin{equation}\label{main_ratio}
\frac{\mathbb{P}(F_1({\cal Q}_1) \cap F_2({\cal Q}_2))}{\mathbb{P}(F_1({\cal Q}_1)) \mathbb{P}(F_2({\cal Q}_2))} = \frac{f(k-2s,t-j_{low},t-j_{up})}{f(k,t,t)},
\end{equation}
where~\(f(a,b,c) := \frac{(a-b)_c}{(a)_c}.\)

Suppose~\(c \geq 1\) and~\(\frac{b+c}{a} < \frac{1}{2}.\) We rewrite
\[f(a,b,c) = \prod_{j=0}^{c-1}\left(1-\frac{b+j}{a}\right) \cdot \prod_{j=1}^{c-1} \left(1-\frac{j}{a}\right)^{-1}\]
and use the estimate~\(e^{-2x} \leq 1-x \leq e^{-x}\) for~\(x< \frac{1}{2}\) to get that
\begin{equation}
f(a,b,c) \leq \prod_{j=0}^{c-1} \exp\left(-\frac{b+j}{a} + \frac{2j}{a}\right) \leq e^{-\frac{bc}{a}+\frac{c^2}{a}} \leq e^{\frac{c^2}{a}}. \label{f_up}
\end{equation}
Similarly
\begin{equation}
f(a,b,c) \geq \prod_{j=1}^{c} \exp\left(-\frac{2b+2j}{a} + \frac{j}{a}\right) = e^{-\frac{2bc}{a}}\exp\left(-\sum_{j=1}^{c}\frac{j}{a}\right) \geq e^{-\frac{2bc}{a} -\frac{c^2}{a}}. \label{f_down}
\end{equation}

We now use~(\ref{f_up}) with~\(a = k-2s,b = t-j_{low}\) and~\(c=t-j_{up}\) to evaluate~\(f(k-2s,t-j_{low},t-j_{up}).\)
We have that
\begin{equation}\label{bca}
\frac{b+c}{a} = \frac{2t-j_1-j_2}{k-2s} \leq \frac{2t}{k-2s} \leq \frac{4t}{k} < \frac{1}{2}
\end{equation}
where the second inequality in~(\ref{bca}) follows from~\(\frac{s}{k}  <\frac{s+t}{k} \leq \frac{1}{4}\)
and the third inequality in~(\ref{bca}) again follows from~\(\frac{t}{k} < \frac{s+t}{k} \leq \frac{1}{8}\) (see statement of Lemma~\ref{pa_lem}).
From~(\ref{f_up}) we therefore have
\begin{equation}
f(k-2s,t-j_{low},t-j_{up}) \leq \exp\left(\frac{(t-j_{up})^2}{k-2s}\right) \leq \exp\left(\frac{t^2}{k-2s}\right)\leq \exp\left(\frac{2t^2}{k}\right) \label{f_num}
\end{equation}
since~\(\frac{s}{k} \leq \frac{s+t}{k} \leq \frac{1}{4}.\) Similarly, using~(\ref{f_down}) with~\(a=k, b=c=t\) we get that
\begin{equation}\label{f_den}
f(k,t,t) \geq \exp\left(-\frac{3t^2}{k}\right)
\end{equation}
Plugging~(\ref{f_num}) and~(\ref{f_den}) into~(\ref{main_ratio}), we get that
\begin{equation}\label{main_ratio2}
\frac{\mathbb{P}(F_1({\cal Q}_1) \cap F_2({\cal Q}_2))}{\mathbb{P}(F_1({\cal Q}_1)) \mathbb{P}(F_2({\cal Q}_2))} \leq \exp\left(\frac{5t^2}{k}\right).
\end{equation}

Substituting~(\ref{main_ratio2}) into~(\ref{ea_split_ax}) and summing over~\({\cal Q}_1\) and~\({\cal Q}_2,\) we get~(\ref{pa_two}).~\(\qed\)

\subsection*{Acknowledgement}
I thank Professors Rahul Roy, C. R. Subramanian and the referee for crucial comments that led to an improvement of the paper. I also thank IMSc and IISER Bhopal for my fellowships.

\bibliographystyle{plain}

\end{document}